\documentclass{gtart}
\usepackage{epsf,latexsym}
\input gtmonout
\volumenumber{1}
\volumeyear{1998}
\volumename{The Epstein birthday schrift}
\pagenumbers{493}{509}
\received{15 November 1997}
\published{29 October 1998}
\papernumber{24}

\def \Z {{\bf Z}}
\def \N {{\bf N}}
\def \R {{\bf R}}
\def \F {{\cal F}}

\def \rtimes {\mbox{$>\!\!\!\lhd$}}
\begin{document}

\title{Hairdressing in groups: a survey of combings\\and formal languages}
\author{Sarah Rees}

\address{University of Newcastle, Newcastle NE1 7RU, UK}
\email{Sarah.Rees@ncl.ac.uk}

\begin{abstract}
A group is combable if it can be represented by a language of words satisfying
a fellow traveller property; an automatic group has a synchronous combing
which is a regular language. 
This article surveys results for combable groups, in particular
in the case where the combing is a formal language.
\end{abstract}

\primaryclass{20F10, 20-04, 68Q40}
\secondaryclass{03D40}
\keywords{Combings, formal languages, fellow travellers, automatic groups}

\maketitle

\cl{\small\it Dedicated to David Epstein on the occasion of his 60th
birthday}

\section{Introduction}
The aim of this article is to survey work generalising the notion of
an automatic group, in particular to classes of groups associated with
various classes of formal languages in the same way that automatic groups
are associated with regular languages.

The family of automatic groups, originally defined by Thurston in an
attempt to abstract certain finiteness properties of the fundamental groups of
hyperbolic manifolds recognised by Cannon in \cite{Cannon}, has been 
of interest for some time.
The defining properties of the family give a geometrical
viewpoint on the groups and facilitate computation
with them; to such a group is associated a set of paths in the Cayley graph
of the group (a `language' for the group) which
both satisfies a geometrical `fellow traveller condition' and, when viewed
as a set of words,
lies in the formal language class of regular languages. (A formal definition is
given in section \ref{automatic}.)
Epstein et al.'s book \cite{ECHLPT} gives a full account; the papers \cite{BGSS} and 
\cite{Farb} are also useful references (in particular, \cite{Farb} is very
readable and non-technical).

The axioms of an automatic group are satisfied by all finite groups, all 
finitely generated free and abelian groups,  word hyperbolic groups,
the fundamental groups of compact Euclidean manifolds, 
and of compact or geometrically finite hyperbolic manifolds 
\cite{ECHLPT, Lang},
Coxeter groups \cite{Brink&Howlett}, braid groups,
many Artin groups \cite{Charney, Charney2, Peifer, Juhasz}, many mapping class
groups \cite{Mosher}, and groups 
satisfying various small cancellation conditions \cite{Gersten&Short}.
However some very interesting groups are not automatic;
the family of automatic groups fails
to contain the fundamental groups of compact 3--manifolds based on
the {\em Nil} or {\em Sol} geometries, and, more generally, fails to contain 
any nilpotent group (probably
also any soluble group) which is not virtually abelian.
This may be surprising since nilpotent groups have very natural
languages, with which computation is very straightforward.

A family of groups which contains the fundamental
groups of all compact, geometrisable 3--manifolds was defined by 
Bridson and Gilman in \cite{Bridson&Gilman}, through a
weakening of both the fellow traveller condition and the
formal language requirement of regularity for automatic groups.
The fellow traveller condition was replaced by an asynchronous
condition of the same type, and the regularity condition by a requirement
that the language be in the wider class of `indexed languages'. 
The class of groups they defined can easily be seen to contain a range of
nilpotent and soluble groups. 

Bridson and Gilman's work suggests that it is sensible to examine other
families of groups, defined in a similar way to automatic groups with respect
to other formal language classes. This paper surveys work on this theme.
It attempts to be self contained,
providing basic definitions and results, but referring the reader
elsewhere for fuller details and proofs.
Automatic groups are defined, and their basic properties described in
section \ref{automatic}; the more general notion of combings is then
explained in section \ref{combings}. A basic introduction to formal languages
is given in section \ref{formal_languages} for the sake of the curious
reader with limited experience in this area. (This section is included to 
set the results of the paper into context, but all or part of it
could easily be omitted on a first reading.)
Section 5 describes
the closure properties of various classes of combable groups,
and section 6 gives examples (and non-examples) of groups with combings
in the classes of regular, context-free, indexed and real-time 
languages.

\rk{Acknowledgment}The author would like to thank the Fakult\"at f\"ur 
Mathematik of the Universit\"at Bielefeld
for its warm hospitality while this work was carried out, and the
Deutscher Akademischer Austauschdienst for financial support.

\section{Automatic groups}
\label{automatic}
Let $G$ be a finitely generated group, and $X$ a finite generating set for $G$,
and define $X^{-1}$ to be the set of inverses of the elements of $X$.
We define a {\em language} for $G$ over $X$ to be a set of {\em words}
over $X$ (that is, products in the free monoid over $X \cup X^{-1}$)
which maps onto $G$ under the natural homomorphism; such a language
is called {\em bijective} if the natural map is bijective.

The group $G$ is automatic if it possesses a language satisfying two
essentially independent conditions, one a geometric `fellow traveller
condition', relating to the Cayley graph $\Gamma$ for $G$ over $X$, the 
other a restriction on the computational complexity of the
language in terms of the formal language class in which the language lives.
Before a precise definition of automaticity can be given, the fellow
traveller condition needs to be explained.

\begin{figure}[htbp]
\begin{center}
\leavevmode
\epsfxsize = 5cm \epsfbox{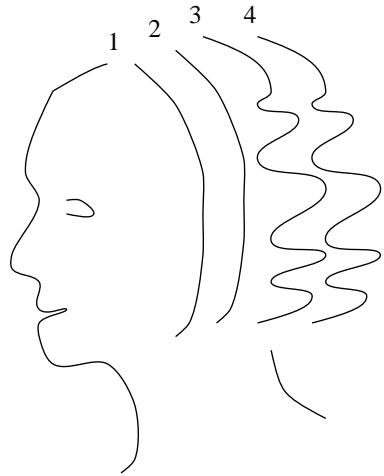}
\caption{ Fellow travellers 
\label{fellow_travellers}}
\end{center}
\end{figure}

Figure \ref{fellow_travellers} gives an informal definition of 
fellow travelling; we give a more formal definition below.
In the figure, the two pairs of
paths labelled 1 and 2, and 3 and 4 synchronously fellow travel at a distance
approximately equal to the length of the woman's nose; the pair of paths
labelled 2 and 3 asynchronously fellow travel at roughly the same distance.
Particles moving at the same speeds along 1 and 2, or along 3 and 4, keep
abreast; but a particle on 3 must move much faster than a particle
on 2 to keep close to it.

More formally let $\Gamma$ be the Cayley graph for $G$ over $X$.
(The vertices of $\Gamma$ correspond to the elements of $G$, and 
an edge labelled
by $x$ leads from $g$ to $gx$, for each $g \in G, x \in X$).
A word $w$ over $X$ is naturally associated with the finite path $\gamma_w$ 
labelled by it and starting at the identity in $\Gamma$.
The path $\gamma_w$ can be parametrised by continuously extending the graph
distance function $d_\Gamma$ (which gives edges length 1); 
where $|w|=d_\Gamma(1,w)$ is the string length of $w$,
for $t\leq |w|$, we define $\gamma_w(t)$ to be a point distance $t$ along
$\gamma_w$ from the identity vertex, and, for $t \geq |w|$, $\gamma_w(t)$ to be
the endpoint of $\gamma_w$. Two paths $\gamma_1$ and $\gamma_2$ of $\Gamma$
are said
to {\em synchronously $K$--fellow travel} if, for all $t\geq 0$, 
$d_\Gamma(\gamma_1(t), \gamma_2(t)) \leq K$, and 
{\em asynchronously $K$--fellow travel} if a strictly increasing positive
valued function 
$h=h_{\gamma_1,\gamma_2}$ can be defined on the positive real numbers,
mapping $[0,l(\gamma_1)+1]$ onto $[0,l(\gamma_2)+1]$,
so that,
for all $t\geq 0$, $d_\Gamma(\gamma_1(t), \gamma_2(h(t))) \leq K$. 

Precisely, $G$ is {\em automatic} if, for some generating set $X$,
$G$ has a language $L$ over $X$ satisfying the following two conditions.
Firstly, for some $K$, and for any $w,v \in L$ for which 
$\gamma_v$ and $\gamma_w$ lead
either to the same vertex or to neighbouring vertices of $\Gamma$,
$\gamma_v$ and $\gamma_w$ synchronously $K$--fellow travel. Secondly
$L$ is regular. A language is defined to be regular if it is the set of words
accepted by a finite state automaton, that is, the most basic form of 
theoretical computer; the reader is referred to section 
\ref{formal_languages} for a crash course on automata theory and 
formal languages.  The regularity of $L$ ensures that
computation with $L$ is easy; the fellow traveller property ensures that the
language behaves well under multiplication by a generator.
Although this is not immediately obvious, the definition of automaticity is
in fact independent of the generating set for $G$; that is,
if $G$ has a regular language over some generating set satisfying the
necessary fellow traveller condition, it has such a language over 
every generating set.

If $G$ is automatic, then $G$ is finitely presented and has quadratic
isoperimetric inequality (that is, for some constant $A$, any loop of
length $n$ in the Cayley graph $\Gamma$ can be divided into at most $An^2$ loops
which are labelled by relators). It follows that $G$ has soluble word problem, 
and in fact there is a straightforward quadratic time algorithm to solve that.

If $G$ is automatic, then so is any subgroup of finite index in $G$, or quotient
of $G$ by a finite normal subgroup, as well as any group in which $G$ is
a subgroup of finite index, or of which $G$ is a quotient by a finite normal
subgroup. The family of automatic groups is also closed under the taking of
direct products, free
products (with finite amalgamation), and HNN extensions (over finite subgroups),
but not under passage to arbitrary subgroups, or under more general products or extensions.

\section{Combings}
\label{combings}

In an attempt to find a family of groups which has many of the good properties
of automatic groups, while also including the examples which are most
clearly missing from that family, we define {\em combable} groups,
using a variant of the first axiom for automatic groups.

Let $G=\langle X \rangle$ be a finitely generated group
with associated Cayley graph $\Gamma$.
We define an {\em asynchronous combing}, or {\em combing} for $G$
to be a language $L$ for $G$ with the property that
for some $K$, and for any $w,v \in L$ for which $\gamma_v$ and $\gamma_w$ lead
either to the same vertex or to neighbouring vertices of $\Gamma$,
$\gamma_v$ and $\gamma_w$ asynchronously $K$--fellow travel; if $G$
has a combing, we say that $G$ is combable. Similarly,
we define a {\em synchronous combing} to be a language for which an analogous
synchronous fellow traveller condition holds; hence automatic groups have 
synchronous combings. Of course, every synchronous combing is also
an asynchronous combing. 

In the above definitions, we have no requirement of
bijectivity, no condition on the length of words in $L$ relative 
to geodesic words, and
no language theoretic restriction. In fact, the term `combing' has been
widely used in the literature, with various different meanings, and some
definitions require some of these properties. Many authors require
combings to be bijective; in \cite{ECHLPT} 
words in the language are required to be quasigeodesic, and in \cite{Gersten}
combings are assumed to be synchronous.

The term `bicombing' is also fairly widely used in the literature,
and so, although we shall not be specifically interested in
bicombability here, we give a definition for the sake of completeness.
Briefly a bicombing is a combing for which words in the language 
related by left multiplication by a generator also satisfy a 
fellow traveller property. Specifically, a combing $L$
is a (synchronous, or asynchronous) {\em bicombing} if paths of the form
$\gamma_v$ and $x\gamma_w$ (synchronously, or asynchronously) fellow travel, 
whenever $\gamma_v,\gamma_w \in L$, 
$x \in X$, and $v=_G xw$, and where $x\gamma_w$ is defined to be the
concatenation of $x$ and a path from $x$ to $xw$ following
edges labelled by the symbols of the word $\gamma_w$.
A group is {\em biautomatic} if it has a synchronous bicombing which is 
a regular language.

Most known examples of combings for non-automatic groups
are not known to be synchronous; certainly this is true of the combings for
the non-automatic groups of compact, geometrisable 3--manifolds found by
Bridson and Gilman. However, in recent and as yet unpublished work,
Bestvina and N. Brady have constructed a synchronous, quasigeodesic
(in fact linear) combing for a non-automatic group.
By contrast, Burillo, in \cite{Burillo}, has shown that none
of the Heisenberg groups 
\begin{eqnarray*}
H_{2n+1}&=&\langle x_1,\ldots x_n,y_1,\ldots y_n,z\,\mid [x_i,y_i]=z,\forall i,\\
& & [x_i,x_j]=[y_i,y_j]=[x_i,y_j]=1,\forall
i,j,i\neq j \rangle
\end{eqnarray*}
or the groups $U_n(\Z)$
of $n$ by $n$ unipotent upper-triangular integer matrices 
can admit synchronous combings by quasigeodesics (all of these groups are
asynchronously combable). Burillo's result was proved by consideration of
higher-dimensional isoperimetric inequalities; the case of $H_3$ had
been previously dealt with in \cite{ECHLPT}.

Let $G$ be a combable group.
Then, by \cite{Bridson} theorem 3.1, $G$ is finitely presented, 
and, by \cite{Bridson} theorems 4.1 and 4.2, $G$ has an exponential
isoperimetric inequality; hence $G$ has soluble word problem
(see \cite{ECHLPT}, theorem 2.2.5).
By \cite{Gersten}, if $G$ has a synchronous, `prefix closed' combing 
(that is, all prefixes of words in the language are in the language),
then $G$ must actually have a quadratic isoperimetric inequality.
Note that, by 
\cite{Kharlampovich} (or see \cite{BGS}),
there are finitely presented class 3 soluble groups which have insoluble
word problem, and so certainly cannot be combable.

For a combing to be of practical use, it must at least be recognisable. 
It is therefore
natural to consider combings which lie in some formal language class, or rather,
which can be defined by some theoretical model of computation.
Automatic groups are associated with the most basic such model, 
that is, with finite state automata and regular languages. In general, 
where $\F$ is a class of formal languages we shall say that a
group is $\F$--combable if it has a combing which is a language in $\F$.
Relevant formal languages are discussed in section \ref{formal_languages}.

An alternative generalisation of automatic groups is discussed in 
\cite{Baumslag&Shapiro&Short}. This approach recognises that the fellow 
traveller condition for a group with language $L$
implies the regularity of the language $L'$ of pairs of words in $L$
which are equal in the group or related by right multiplication by a generator,
and examines what happens when both $L$ and $L'$ are allowed to lie
in a wider language class (in this particular case languages are considered
which are intersections of context-free languages, and hence
defined by series of pushdown automata). Some of the consequences of such a
generalisation are quite different from those of the case of combings; 
for example, such groups need not be
finitely presented.

\section{Hierarchy of computational machines and formal languages}
\label{formal_languages}
Let $A$ be a finite set of symbols, which we shall call an {\em alphabet}.
We define a {\em language} $L$ over $A$ to be a set of finite strings (words)
over $A$, that is a subset of $A^* = \cup_{i\in \N}A^i$.
We define a {\em computational machine} $M$ for $L$ to be a device 
which can be used to recognise the words in $L$, as follows. 
Words $w$ over $A$ can be input to $M$ one at a time for processing.  
If $w$ is in $L$, then the processing of $w$ terminates after some finite time,
and $M$ identifies $w$ as being in $L$; if $w$ is not in $L$, then either $M$
recognises this after some time, or $M$ continues processing $w$ indefinitely.
We define $L$ to be a {\em formal language} if it can be
recognised by a computational machine; machines of varying complexity define
various families of formal languages. 

We shall consider various different types of computational machines. Each
one can be described in terms of two basic components, namely a finite set $S$
of {\em states}, between which $M$ fluctuates,
and (for all but the simplest machines) a possibly infinite {\em memory} 
mechanism. 
Of the states of $S$, one is identified as a {\em start state} and some are
identified as {\em accept states}. 
Initially (that is, before
a word is read) $M$ is always in the start state; 
the accept states are used by $M$ to help it in its decision process, possibly
(depending on the type of the machine)
in conjunction with information retrieved from the memory.

We illustrate the above description with a couple of examples of formal
languages over the alphabet $A = \{-1,1\}$, and machines which recognise them. 

We define $L_1$ to be the language over $A$
consisting of all strings containing an even number of $1$'s. This language
is recognised by a very simple machine $M_1$ with two states and no additional
memory. $S$ is the set $\{even,\ odd \}$; $even$ is the start state
and only accept state.
$M_1$ reads each word $w$ from left to right, and switches state
each time a $1$ is read. The word $w$ is accepted if $M_1$ is in the
state $even$ when it finishes reading $w$. $M_1$ is an example
of a (deterministic) finite state automaton.

We define $L_2$ to be the language over $A$ consisting of all strings 
containing an equal number of $1$'s and $-1$'s. This language
is recognised by a machine $M_2$ which reads an input word $w$ from left to 
right,
and keeps a record at each stage of the sum of the digits so far read; $w$
is accepted if when the machine finishes reading $w$ this sum is equal to
$0$. For this machine the memory is the crucial component (or rather, 
the start state is the only state). The language $L_2$ cannot be recognised 
by a machine without
memory. $M_2$ is an example of a pushdown automaton.

A range of machines and formal language families, ranging from the simplest
finite state automata and associated regular (sometimes known as rational)
languages to the Turing machines and recursively enumerable languages,
is described in \cite{Hopcroft&Ullman}; a treatment
directed towards geometrical group theorists is provided by \cite{Gilman}.
One-way nested stack automata and real-time Turing machines (associated
with indexed languages and real-time languages respectively) are also of 
interest to us in this article, and are discussed in \cite{Aho,Aho2} and
in \cite{Rabin,Rosenberg}. We refer the reader to those papers for details,
but below we try to give an informal overview of relevant machines and
formal languages. 

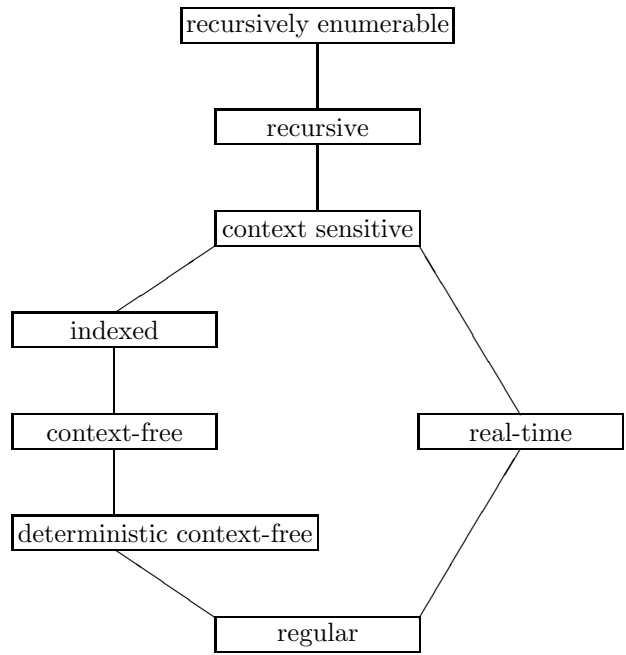
\begin{figure}[htb]
\centering
\bigskip
\linethickness{0.4pt}
\unitlength 0.9mm
\begin{picture}(90.00,100.00)\small
\put(30.00,00.00){\framebox(30.00,05.00){regular}}
\put(30.00,05.00){\line(-3,2){15.00}}
\put(60.00,05.00){\line(3,5){15.00}}
\put(00.00,15.00){\framebox(45.00,05.00){deterministic context-free}}
\put(15.00,20.00){\line(0,1){10.00}}
\put(00.00,30.00){\framebox(30.00,05.00){context-free}}
\put(15.00,35.00){\line(0,1){10.00}}
\put(00.00,45.00){\framebox(30.00,05.00){indexed}}
\put(15.00,50.00){\line(3,2){15.00}}
\put(60.00,30.00){\framebox(30.00,05.00){real-time}}
\put(75.00,35.00){\line(-3,5){15.00}}
\put(30.00,60.00){\framebox(30.00,05.00){context sensitive}}
\put(45.00,65.00){\line(0,1){10.00}}
\put(30.00,75.00){\framebox(30.00,05.00){recursive}}
\put(45.00,80.00){\line(0,1){10.00}}
\put(25.00,90.00){\framebox(40.00,05.00){recursively enumerable}}
\end{picture}
\caption{Inclusions between formal language classes 
\label{hierarchy}}
\end{figure}

Figure \ref{hierarchy} shows known inclusions between the formal language
classes which we shall describe.

We continue with descriptions of various formal language classes; these
might be passed over on a first reading.

\subsection{Finite state automata and regular languages}
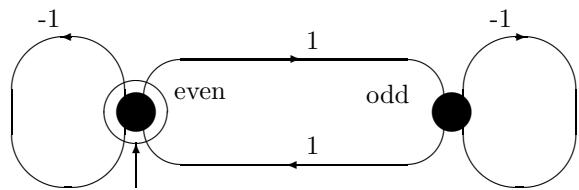
\begin{figure}[htb]
\centering
\bigskip
\linethickness{0.4pt}
\unitlength 1mm
\begin{picture}(90.00,30.00)\small
\put(40.00,22.00){\oval(40.00,10.00)[t]}
\put(40.00,27.00){\vector(1,0){1}}
\put(40.00,27.00){\makebox(05.00,05.00){1}}
\put(40.00,18.00){\oval(40.00,10.00)[b]}
\put(40.00,13.00){\vector(-1,0){1}}
\put(40.00,13.00){\makebox(05.00,05.00){1}}
\put(10.00,20.00){\oval(15.00,20.00)}
\put(25.00,20.00){\makebox(05.00,05.00){even}}
\put(10.00,30.00){\vector(-1,0){1}}
\put(05.00,30.00){\makebox(05.00,05.00){-1}}
\put(19.00,20.00){\circle{8}}
\put(19.00,20.00){\circle*{6}}
\put(19.00,10.00){\vector(0,1){6}}
\put(70.00,20.00){\oval(15.00,20.00)}
\put(50.00,20.00){\makebox(05.00,05.00){odd}}
\put(70.00,30.00){\vector(1,0){1}}
\put(65.00,30.00){\makebox(05.00,05.00){-1}}
\put(61.00,20.00){\circle*{6}}
\end{picture}
\vglue -9mm
\caption{The finite state automaton $M_1$ 
\label{automaton}}
\end{figure}

A set of words over a finite alphabet is defined to be a {\em regular}
language precisely if it is the language defined by a finite state automaton.
A {\em finite state automaton} is a machine without memory,
which moves through the states of $S$ as it reads words over $A$
from left to right. 
The simplest examples are the so-called {\em deterministic} finite
state automata.
For these a transition function 
$\tau\co  S \times A \rightarrow S$ determines passage between states;
a word $w=a_1\ldots a_n$ ($a_i \in A$) is accepted
if for some sequence of states $s_1,\ldots s_n$, of which $s_n$ is an accept 
state, for each $i$, $\tau(s_{i-1},a_i) = s_i$.
Such a machine is probably best understood when viewed as a finite,
directed, edge-labelled graph (possibly with loops and multiple edges), 
of which the states are vertices.
The transition $\tau(s,a) = s'$ is then represented by an edge labelled by
$a$ from the vertex $s$ to the vertex $s'$.
At most one edge with any particular label
leads from any given vertex (but since dead-end non-accept states can easily
be ignored, there may be less that $|A|$ edges out of a vertex, and further,
several edges with distinct labels might connect the same pair of vertices).
A word $w$ is accepted if
it labels a path through the graph from the start vertex/state $s_0$
to a vertex which is marked as an accept state.
Figure \ref{automaton} gives such a graphical description for the machine
$M_1$ described at the beginning of section \ref{formal_languages}.
In such a figure,
it is customary to ring the vertices which represent accept states,
and to point at the start state with a free arrow, hence the state
$even$ is recognisable in this figure as the start state and sole accept state.

A {\em non-deterministic} finite state automaton is
defined in the same way as a deterministic finite state automaton
except that the transition function $\tau$
is allowed to be multivalued. A word $w$ is accepted if some (but
not necessarily all) sequence of transitions following the symbols of
$w$ leads to an accept
state. The graphical representation of a non-deterministic machine
may have any finite number of edges with a given label from
each vertex. In addition, further edges
labelled by a special symbol $\epsilon$ may allow the machine to leap, without
reading from the input string, from one state to another, in a so-called
$\epsilon$--move.

Given any finite state automaton, possibly with multiple edges from
a vertex with the same label, possible with $\epsilon$--edges,
a finite state automaton defining the same language can be constructed in which
neither of these possibilities occur. Hence, at the level of finite state
automata, there is no distinction between the deterministic and
non-deterministic models. However, for other classes of machines (such
as for pushdown automata, described below) non-determinism increases the power
of a machine.

\subsection{Turing machines and recursively enumerable languages}
The {\em Turing machines}, associated with the {\em recursively enumerable} 
languages,
lie at the other end of the computational spectrum from finite state automata,
and are accepted as providing a formal definition of computability.
In one of the simplest models (there are many equivalent models)
of a Turing machine, we consider the input word to be written on a section of
a doubly-infinite tape, which is read through a movable {\em tape-head}.
The tape also serves as a memory device.
Initially the tape contains only the input word $w$, the tape-head
points at the left hand symbol of that word, and the machine is in the
start state $s_0$. Subsequently, the tape-head may move both right and
left along the tape (which remains stationary). At any stage, the tape-head
either reads the symbol from the section of tape at which it currently points
or observes that no symbol is
written there. Depending on the state it is currently in, and what it
observes on the tape, the machine changes state, writes a new symbol
(possibly from $A$, but possibly one of finitely many other
symbols, or blank) onto the tape, and either halts, or moves its tape-head 
right or left one position. 
The input word $w$ is accepted if the machine eventually halts in an
accept state; it is possible that the machine may not halt on all input.

Non-deterministic models, where the machine may have a choice of moves
in some situations (and accepts a word if some allowable sequence
of moves from the obvious initial situation leads it to halt in an accept 
state), and models with any finite number of extra tapes and tape-heads, 
are all seen to be 
equivalent to the above description, in the sense that they also define 
the recursively enumerable languages.

\subsection{Halting Turing machine and recursive languages}
A {\em halting Turing machine} is a Turing machine which halts on all input; thus
both the language of the machine and its complement are recursively 
enumerable. A language
accepted by such a machine is defined to be a {\em recursive language}.

\subsection{Linear bounded automaton and context sensitive languages}
A {\em linear bounded automaton} is a non-deterministic Turing machine
whose tape-head is only allowed to move through the piece of tape which 
initally contains the input word; special symbols, which cannot be overwritten,
mark the two ends of the tape.  Equivalently (and hence the name), 
the machine is restricted to a piece of tape whose length is a 
linear function of the length of the input word.
A language accepted by such a machine is defined to be a {\em context
sensitive language}.

\subsection{Real-time Turing machines and real-time languages}
A {\em real-time Turing machine} is most easily described as a
deterministic Turing machine with
any finite number of doubly-infinite tapes (one of which initially contains
the input, and the others of which are initially empty), which halts
as it finishes reading its input. Hence such a machine processes its input in
`real time'. 

A `move' for this machine consists of an operation of each of the tape heads,
together with a state change, as follows.
On the input tape, the tape-head reads the symbol to which it currently points, 
and then moves one place to the right.
On any other tape, the tape-head reads the symbol (if any) to which it currently
points, prints a new symbol (or nothing), and then either moves right,
or left, or stays still.
The machine changes to a new state, which depends on its current state, 
and the symbols read from the tapes.  
When the tape-head on the input head has read the last symbol of the input,
the whole machine halts, and the input word is accepted if the machine
is in an accept state.

A language accepted by such a machine is defined to be a
{\em real-time language}.
$\{a^nb^nc^n:n \in \N\}$ is an example \cite{Rosenberg}.
Examples are descibed in \cite{Rosenberg}
both of real-times languages which do not lie in the
class of context-free languages (described below), 
and of (even deterministic) context-free 
languages which are not real-time.

\subsection{Pushdown automata and context-free languages}
A {\em pushdown automaton} can be described as a Turing machine with a 
particularly restricted operation on its tape, but it is probably easier to
visualise as a machine formed by adding an infinite stack 
(commonly viewed as a spring-loaded pile of plates in a canteen)
to a (possibly non-deterministic) finite state automaton. 
Initially the stack contains a single start symbol.
Only the top symbol of the stack can be accessed at
any time, and information can only be appended to the top of the stack.
The input word $w$ is read from left to right. 
During each move, the top symbol of the stack is removed from the stack,
and a symbol from $w$ may be read, or may not.
Based on
the symbols read, and the current state of the machine, the machine moves
into a new state,
and a string of symbols (possibly empty) from a finite 
alphabet is appended to the top of the stack. The word $w$ is accepted if after 
reading it the machine may be in an accept state. 
The language accepted by a pushdown automaton is defined to be a 
{\em context-free language}. 

The machine $M_2$ described towards the beginning of this section
can be seen to be a pushdown automaton as follows. The `sum so far'
is held in memory as either a sequence of $+1$'s or as a sequence
of $-1$'s with the appropropriate sum. When the top symbol on the
stack is $+1$ and a $-1$ is read from the input tape, 
the top stack symbol is removed,
and nothing is added to the stack.  When the top symbol on the
stack is $-1$ and a $+1$ is read from the input tape, 
the top stack symbol is removed,
and nothing is added to the stack. Otherwise, the top stack
symbol is replaced, and then the input symbol is
added to the stack. Hence the language $L_2$ recognised by $M_2$
is seen to be context-free. Similarly so is the language
$\{a^nb^n: n \in N\}$ over the alphabet $\{a,b\}$.
Neither language is regular. For symbols $a,b,c$, the language
$\{a^nb^nc^n: n \in N\}$ is not context-free.

A pushdown automaton is deterministic if each input word $w$ defines
a unique sequence of moves through the machine. This does not in fact mean
that a symbol of $w$ must be read on each move, but rather that the decision
to read a symbol from $w$ at any stage is determined by the symbol read from
the stack and the current state of the machine.
The class of deterministic context-free languages
forms a proper subclass of the class of context-free languages,
which contains both the examples of context-free languages given above.
The language consisting of all words of the form $ww^R$ over some
alphabet $A$ (where $w^R$ is the reverse of $w$) is 
non-deterministic context-free \cite{Hopcroft&Ullman}, but is
not deterministic context-free.

\subsection{One-way nested stack automata and indexed languages}
A {\em one-way nested stack automaton} is probably most easily viewed
as a generalisation of a pushdown automaton, that is, as a non-deterministic
finite state automaton with an attached nest of stacks, rather than a single
stack.
The input word is read from left to right (as implied by the term `one-way').
In contrast to a pushdown automaton,
the read/write tape-head of this machine is allowed some movement through
the system of stacks.
At any point of any stack to which the tape-head has access it can read, and
a new nested stack can be created; while at the top of any stack it can also 
write,
and delete. 
The tape-head can move down through any stack, but its
upward movement is  restricted; basically it is not allowed to move
upwards out of a non-empty stack.

The language accepted by a one-way nested stack automaton is defined to be an
{\em indexed language}. 
For symbols $a,b,c$,
the languages  $\{a^nb^nc^n:n \in \N\}$,
$\{ a^{n^2}: n\geq 1\},\{a^{2^n}:n\geq 1\}$ and
$\{ a^nb^{n^2}: n \geq 1\}$ are indexed \cite{Hopcroft&Ullman}, but
$\{ a^{n!} : n \geq 1\}$ is not \cite{Hayashi}, nor is
$\{ (ab^n)^n: n \geq 1 \}$ \cite{Gilman2, Hayashi}.

\section{From one $\F$--combing to another}
\label{closure}
Many of the closure properties of the family of automatic
groups also hold for other classes of combable groups, often for
synchronous as well as asynchronous combings.

In the list below we assume that $\F$ is either the set of all
languages over a finite alphabet, or is one of the classes of formal 
languages described in section \ref{formal_languages}, that is that 
$\F$ is one of the
regular languages, context-free languages, indexed languages, 
context-sensitive languages, real-time languages, recursive languages, 
or recursively enumerable languages. 
(These results for all but real-time languages are proved in \cite{Bridson&Gilman} and \cite{Rees2}, and for real-time languages in \cite{GHR}.)
Then just as for automatic groups, we have all the following results:
\begin{itemize}
\item If $G$ has a synchronous or asynchronous $\F$--combing
then it has such a combing over any generating set.
\item Where $N$ is a finite, normal subgroup of $G$, 
and $G$ is finitely generated, 
then $G$ is synchronously or asynchronously $\F$--combable if and only
the same is true of $G/N$.
\item Where $J$ is a finite index subgroup of $G$, 
then $G$ is synchronously or asynchronously $\F$--combable if and only
if the same is true of $J$.
\item If $G$ and $H$ are both asynchronously $\F$--combable then so
are both $G\times H$ and $G \ast H$.
\end{itemize}

A crucial step in the construction of combings for 3--manifold groups
in \cite{Bridson&Gilman} is a construction of Bridson in  \cite{Bridson2};
combings for $N$ and $H$ can be put together to give an asynchronous
combing for a split-extension of the form $N \rtimes H$
provided that $N$ has a combing which is particularly stable under
the action of $H$. 
The set of all geodesics in a word hyperbolic group has that stability,
and is a regular language;
hence,  for any of the language classes $\F$ considered in this section,
any split extension of a word hyperbolic group by an $\F$--combable
group is $\F$--combable. The free abelian group $\Z^n$ also possesses
a combing with the necessary stability;  hence all split extensions of
$\Z^n$ by combable groups are asynchronously combable. It remains only
to ask in which language class these combings lie.

Stable combings for $\Z^n$ are constructed by Bridson in \cite{Bridson2}
as follows. $\Z^n$ is seen embedded as a lattice in $\R^n$, and the
group element $g$ is then represented by
a word which, as a path through the
lattice, lies closest to the real line joining the point $0$ to the
point representing $g$. For some group elements there is a selection of
such paths; a systematic choice can clearly be made.
It was proved in \cite{Bridson&Gilman} that $\Z^2$ has a combing of this type
which is an indexed language; hence all split extensions of the
form $\Z^2 \rtimes \Z$ were seen to be indexed combable. It followed from this
that the fundamental groups of all compact, geometrisable 3--manifolds
were indexed combable;
for these are all commensurable with free products 
of groups which are either automatic or finite extensions of $\Z^2 \rtimes \Z$.

It is unclear whether or not the corresponding combing for $\Z^n$ is also
an indexed language when $n>2$. Certainly  it is a real-time language
\cite{GHR}. Hence many split extensions of the form
$\Z^n \rtimes H$ are seen to have asynchronous combings
which are real-time languages. We give some examples in the
final section.

\section{Combing up the language hierarchy}
\subsection{Regular languages}
A group with a synchronous regular combing is, by definition, automatic.
More generally, a group with a regular combing is called {\em asynchronously
automatic} \cite{ECHLPT}. It is proved in \cite{ECHLPT} that the 
asynchronicity of an asynchronously automatic group is bounded; that is the
relative speed at which particles must move along two fellow-travelling
words in order to keep apace can be kept within bounds.
The Baumslag--Solitar groups 
\[ G_{p,q} = \langle a, b \mid ba^p = a^q b \rangle \]
are asynchronously automatic, but not automatic, for $p \neq \pm q$ 
(see \cite{ECHLPT, Rees}), and automatic for $p = \pm q$.

It is proved in \cite{ECHLPT} that a nilpotent group which is not 
abelian-by-finite cannot be asynchronously automatic. From this it follows
that the fundamental groups of compact manifolds based on the $Nil$
geometry cannot be asynchronously automatic; N. Brady proved
that the same is true of groups of the compact manifolds based on
the $Sol$ geometry \cite{Brady}.

\subsection{Context-free languages}
No examples are currently known of non-automatic groups with
context-free combings.
It is proved in \cite{Bridson&Gilman} that a nilpotent group which is not
abelian-by-finite cannot have a bijective 
context-free combing; however it remains open whether a context-free
combing with more that one representative for some group elements
might be possible.
\subsection{Indexed languages}
Bridson and Gilman proved that the fundamental group
of every compact geometrisable 3--manifold (or orbifold) is
indexed combable. By the results of \cite{Brady, ECHLPT, 
Bridson&Gilman} described above for regular and context-free 
combings, this result must be close to being best possible.

It follows immediately from Bridson and Gilman's results that
a split extension of $\Z^2$ by an indexed combable (and so, certainly
by an automatic) group is again indexed combable.

\subsection{Real-time languages}
Since the stable combing of $\R^n$ described in section \ref{closure}
is a real-time language \cite{GHR}, it follows that any split
extension over $\Z^n$ of a real-time combable group is real-time combable.
Hence (see \cite{GHR}),
any finitely generated class 2 nilpotent group with cyclic 
commutator subgroup is real-time combable, and also any 3--generated
class 2 nilpotent group. Further the free class 2 nilpotent
groups, with presentation,
\[ \langle x_1,\ldots x_k \mid [[x_i,x_j],x_k],\,\forall i,j,k \rangle,\]
as well as the $n$--dimensional Heisenberg groups
and the groups of $n$--dimen\-sional, unipotent upper-triangular 
integer matrices, can all be expressed as split extensions over
free abelian groups, and hence are real-time combable.
It follows that
any polycyclic-by-finite group (and so,
in particular, any finitely generated nilpotent group) embeds
as a subgroup in a real-time combable group.

Torsion-free polycyclic metabelian groups with centre disjoint 
from their commutator subgroup are far from being nilpotent, but are
also real-time combable (see \cite{GHR}). Such groups split 
over their commutator subgroup, by a theorem of \cite{Robinson}.
An example is provided by the group
\[ \langle x,y,z \mid  yz = zy, y^x = yz,  z^x = y^2z \rangle \]
which is certainly not automatic (it has exponential isoperimetric
inequality). In fact this group is also indexed combable, since
it is of the form $\Z^2 \rtimes \Z$.

\Addresses\recd


\begin{thebibliography}

\itemsep 1.5pt

\bibitem{Aho} {\bf Alfred V Aho}, {\it Indexed grammars -- an extension
of context-free grammars}, J. Assoc. Comp. Mach. 15 (1968) 647--671

\bibitem{Aho2} {\bf Alfred V Aho}, {\it Nested stack automata},
J. Assoc. Comp. Mach. 16 (1969) 383--406

\bibitem{BGSS} {\bf G Baumslag}, {\bf S\,M Gersten}, {\bf M Shapiro},
{\bf H Short}, {\it Automatic groups and amalgams}, Journal of
Pure and Applied Algebra 76 (1991) 229--316

\bibitem{BGS} {\bf G Baumslag}, {\bf D Gildenhuys}, {\bf R Strebel},
{\it Algorithmically insoluble problems about finitely presented
solvable groups, Lie and associative algebras. I}, Journal of
Pure and Applied Algebra 39 (1986) 53--94


\bibitem{Baumslag&Shapiro&Short} {\bf Gilbert Baumslag}, {\bf Michael
Shapiro}, {\bf Hamish Short}, {\it Parallel poly pushdown groups},
Jorunal of Pure and Applied Algebra, to appear

\bibitem{Brady} {\bf Noel Brady}, {\it The geometry of asynchronous automatic 
structures on groups}, PhD thesis, UC Berkeley (1993)

\bibitem{Bridson} {\bf Martin R Bridson},
{\it On the geometry of normal forms in discrete groups},
Proc. London Math. Soc. (3) 67 (1993) 596--616

\bibitem{Bridson2} {\bf Martin R Bridson},
{\it Combings of semidirect products and $3$--manifold groups},
Geometric and Functional Analysis 3 (1993) 263--278

\bibitem{Bridson&Gilman} {\bf M\,R Bridson}, {\bf R\,H Gilman}, 
{\it Formal language theory and the geometry of $3$--manifolds},
Commentarii Math. Helv. 71 (1996) 525--555

\bibitem{Brink&Howlett} {\bf Brigitte Brink}, {\bf Robert Howlett},
{\it A finiteness  property of Coxeter groups}, Math. Ann. 296
(1993) 179--190

\bibitem{Burillo} {\bf Jos\'e Burillo},   
{\it Lower bounds of isoperimetric functions for nilpotent groups},
DIMACS Ser. Discrete Math. Theoret. Comput. Sci. (25) 1--8

\bibitem{Cannon} {\bf J\,W Cannon}, {\it The combinatorial structure
of cocompact discrete hyperbolic groups}, Geom. Dedicata 16 (1984) 123--148

\bibitem{Charney} {\bf Ruth Charney}, {\it Artin groups of finite type
are biautomatic}, Math. Ann. 292 (1992) 671--683

\bibitem{Charney2} {\bf Ruth Charney},
{\it Geodesic automation and growth functions of Artin groups}, Math. Ann. 301 (1995) 307--324

\bibitem{ECHLPT} {\bf David B\,A Epstein}, {\bf J\,W Cannon}, {\bf
D\,F Holt}, {\bf S Levy}, {\bf M\,S Patterson}, {\bf W Thurston}, {\it
Word processing in groups}, Jones and Bartlett, (1992)

\bibitem{Farb} {\bf Benson Farb}, {\it Automatic groups: a guided
tour}, Enseign. Math. (2) 38 (1992) 291--313

\bibitem{Gersten} {\bf S\,M Gersten}, {\it Bounded cocycles and
combings of groups}, Internat. J. Algebra Comput. 2 (1992) 307--326

\bibitem{Gersten&Short} {\bf S\,M Gersten}, {\bf H\,B Short}, {\it Small
cancellation theory and automatic groups}, Inv. Math. 102 (1990) 305--334

\bibitem{Gilman} {\bf Robert H Gilman}, {\it Formal languages and infinite
groups}, DIMACS Ser. Discrete Math. Theoret. Comput. Sci. 25 AMS
(1996) 27--51

\bibitem{Gilman2} {\bf Robert H Gilman}, {\it A shrinking lemma for
indexed languages}, Theoret. Comput. Sci. 163 (1996) 277--281

\bibitem{GHR} {\bf Robert H Gilman}, {\bf Derek F Holt}, {\bf Sarah
Rees}, {\it Combing nilpotent and polycylic groups},
Internat. J. Algebra Comput. to appear

\bibitem{Hayashi} {\bf T Hayashi}, {\it On derivation trees of indexed
grammars}, Publ. RIMS Kyoto Univ. 9 (1973) 61--92

\bibitem{Hopcroft&Ullman} {\bf John E Hopcroft}, {\bf Jeffrey D
Ullman}, {\it Introduction to automata theory, languages and computation}, Addison--Wesley, (1979)

\bibitem{Juhasz} {\bf Arye Juh\'asz}, {\it Large Artin groups are biautomatic}, preprint

\bibitem{Kharlampovich} {\bf O\,G Kharlampovich}, {\it A finitely
presented group with unsolvable word problem (in Russian)},
Izv. Akad. Nauk SSSR Ser. Mat. 45 (1981) 852--873

\bibitem{Lang} {\bf U Lang}, {\it Quasigeodesics outside horoballs}, 
Geom. Dedicata 63 (1996) 205--215

\bibitem{Mosher} {\bf Lee Mosher}, {\it Mapping class groups are automatic},
Annals of Math. 142 (1995) 303--384

\bibitem{Peifer} {\bf David Peifer}, {\it Artin groups of extra-large type are
automatic}, J. Pure Appl. Alg. 110 (1996) 15--56

\bibitem{Rabin} {\bf Michael O Rabin}, {\it Real time computation},
Israel J. Math. 1 (1963) 203--211

\bibitem{Rees} {\bf Sarah Rees}, {\it Automatic groups associated with
word orders other than shortlex}, Int. J. Alg. Comp. (to appear)

\bibitem{Rees2} {\bf Sarah Rees}, {\it A language theoretic analysis of combings}, preprint

\bibitem{Robinson} {\bf Derek J\,S Robinson}, {\it Splitting theorems
for infinite groups}, Symposia Mathematica 17, Convegni del Novembre e
del Dicembre 1973, Academic Press, (1976) 441--470

\bibitem{Rosenberg} {\bf A Rosenberg}, {\it Real-Time Definable
Languages}, J. Assoc. Comput. Mach. 14 (1967) 645-662


\end{thebibliography}
\end{document}